\documentclass[12pt]{article}

\usepackage{amsfonts}
\usepackage{amsmath}
\usepackage{amssymb}
\usepackage{amscd}
\usepackage{verbatim}
\usepackage{latexsym}

\newcommand{\A}{\mathbb{A}}
\newcommand{\Q}{\mathbb{Q}}

\newcommand{\aQ}{\overline{\mathbb{Q}}}
\newcommand{\Z}{\mathbb{Z}}

\newcommand{\aFp}{\overline{\mathbb{F}_{p}}}

\newcommand{\aQp}{\overline{\mathbb{Q}_{p}}}
\newcommand{\Qp}{\mathbb{Q}_{p}}
\newcommand{\GL}{\mathrm{GL}}

\newcommand{\rhobar}{\overline{\rho}}
\newtheorem{theorem}{Theorem}[section]
\newtheorem{lemma}[theorem]{Lemma}
\newtheorem{prop}[theorem]{Proposition}

\newtheorem{cor}[theorem]{Corollary}

\newenvironment{proof}{\smallskip\noindent{\bf Proof.}\
}{$\Box$\smallskip\goodbreak}

\newcommand{\FF}{ {\mathbb F} }

\def\rhobar{ {\bar {\rho} } }
\def\rhob{ {\bar {\rho} } }

\newcommand{\Galois}{\mathrm{Gal}}
\newcommand{\Gal}{\Galois( \bar{ {\mathbb Q}}/{\mathbb Q})}
\newcommand{\F}{{\mathbb F}}

\begin{document}

\title{On Serre's reciprocity conjecture for 2-dimensional mod $p$
 representations of $\Gal$ \\ }

\author{Chandrashekhar Khare, University of Utah and
TIFR \thanks{partially supported
by NSF Grant DMS - 0355528} , \\
Jean-Pierre Wintenberger, Universit\'e Louis Pasteur
}


\maketitle

\begin{abstract}
We first prove  the existence of {\it minimally
  ramified}  $p$-adic lifts of 
2-dimensional mod $p$  representations $\rhobar$, that are odd and
irreducible, of the absolute Galois group of $\Q$,  in many cases. This is
predicted by Serre's conjecture that such representations arise from
newforms of optimal level and weight.

Using these minimal lifts, and arguments using compatible systems,
we  prove some cases of Serre's conjectures in low levels and
weights. For instance we prove that there are no irreducible $(p,p)$ type
group schemes over the rational integers. We prove that a $\rhobar$ as above
of Artin conductor 1 and Serre weight 12 arises from the Ramanujan Delta-function.

In the last part of the paper we present arguments that reduce Serre's
conjecture on
modularity of representations $\rhobar$ as above to proving modularity
lifting theorems of the type pioneered
by Wiles. While these modularity lifting results are not known as yet they
might be relatively accessible because of the basic method of Wiles and
Taylor-Wiles and recent developments in the p-adic Langlands programme
initiated by Breuil.

\end{abstract}

\tableofcontents

\newcommand{\vep}{{\varepsilon}}

\section{Introduction}

Fix an irreducible 2 dimensional odd mod $p$ representation
$\rhob:\Gal \rightarrow GL_2(\FF)$ with $\FF$ a finite field of
characteristic
$p$ which we assume is odd. We say that such a representation is of {\it
Serre
type}. Serre in
\cite{Serre3} conjectures that all such
representations arise from newforms. The main technique which is presented
in this paper results in the reduction of proving this conjecture in many
cases
to proving generalisations of modularity lifting theorems of the type
pioneered by Wiles. By modularity lifting theorems we mean proving
that certain characteristic 0 lifts $\rho$ of modular $\rhobar$, i.e.,
$\rhobar$
such that $\rhobar^{ss}$ arises from a modular form (so by convention if
$\rhobar$ is reducible it is modular), are again modular.
The ideas of this paper also lead to unconditional proofs of
Serre's conjecture
in low levels and low weights. (Below we say a
Serre type $\rhob$ is of level one if it is unramified outside $p$).

The 3 main results of this paper are:

\begin{enumerate}

\item Liftings of $\rhobar$ to {\it minimally ramified representations} (see
 Theorem \ref{existminimallift}).

\item Proof of Serre's conjectures in low levels and weights (see Theorem
\ref{lowlevel}, \ref{lowweight} and \ref{weight12}).

\item The reduction of many cases of Serre's conjecture to modularity
lifting
theorems. More specifically:

- ($p$ bigger than $3$) The reduction of Serre's conjecture
for representations  with odd, prime to 3, squarefree Artin conductor and
Serre
weight 2
to the level one case, assuming modularity lifting theorems
for semistable representations between weights 2 and $p+1$
(see Theorem \ref{levelone}).

- The reduction of Serre's conjecture in the level one case to modularity
lifting theorems (see Theorem \ref{bertrand}).
\end{enumerate}

These reductions are curious, and at first sight may seem surprising,
because
modularity lifting results, as their name
suggests,  {\it assume} residual modularity.
So in fact we show that modularity lifting results in principle (and when
proved
in sufficient generality even actually) imply residual modularity. For this
we need to know modularity lifting results a little beyond what is
known
presently. (We need modularity of lifts of modular, possibly reducible or
dihedral,
odd 2 dimensional mod $p$ representations when the lifts are either
crystalline and of weight between 2 and $2p$, or semistable of weight
between 2
and $p+1$).
One might hope, not too extravagantly, to prove these extensions
using the basic method of Wiles and Taylor-Wiles,
especially because of the work of Breuil, and related developments
in the $p$-adic Langlands program due to Berger, Li, Zhu and Mezard
(see \cite{Breuil},
\cite{BM},
\cite{BLZ}), and  recent developments of the $R=T$ machinery
due to Kisin (see \cite{Kisin}).

\vspace{.5cm}

The main techniques  of this paper are:

- The use of known modularity lifting techniques pioneered by
Wiles, and Taylor and Wiles
and crucial
further developments of Taylor, see \cite{[T02]} and \cite{[T03]},
which prove a potential version of Serre's conjecture (and that
introduced  the  technique of reducing this to modularity lifting results).
This leads to proving finiteness of minimal deformation rings
attached to $\rhob:\Gal \rightarrow GL_2(\FF)$ (to be defined below:
especially
what we mean by this at $p$ is delicate and context-dependent). For this,
we use
base change arguments of the type
used by de Jong in \cite{[dJ01]} and Part II of \cite{BoKh}, and then
results of B\"ockle in the
appendix to \cite{[K03]}.
This gives that these minimal deformation rings
are finite, flat complete intersections. This leads to
existence of a minimal $p$-adic lift $\rho$ of $\rhob$.
This argument was known in principle
to the authors of \cite{[KR03]}.
We emphasise that this is an existence theorem in a sense, and to us it
seems
very hard to produce minimal lifts {\it by hand}.

- Then arguments of Dieulefait and the second author, see
 \cite{D}
and \cite{Wint} and also \cite{[T03]},
are used to make
$\rho$ part of a {\it minimal strictly compatible system}.
Then it is immediate to
prove that there are no irreducible $(p,p)$ finite flat group schemes over
$\Z$. Refining these ideas and using refinements of Fontaine's results
 in
\cite{Fontaine} as
in work of Brumer and Kramer \cite{BK}, and Schoof \cite{Schoof1}, leads to
proof of Serre's conjecture in many low level and low weight cases (
for example Serre type $\rhob$ of Serre weight 12, and level 1, always
arises from the Ramanujan $\Delta$ function
which uses a beautiful result of Schoof in \cite{Schoof1}).

- Arguments using existence of minimal compatible systems of different types
to reduce in many cases the semistable case of Serre's conjecture to
the level 1 case (this step can be called using {\it switching} via
compatible
systems to {\it kill ramification}) assuming modularity lifting results
(see Theorem \ref{levelone}). Another use of these arguments, together
with
the prime number theorem, in fact a much weaker version of it known as {\it
Bertrand's postulate}, produces the last reduction mentioned above (see
Theorem \ref{bertrand}).

\vspace{.5cm}

It is to be noted that Ramakrishna in \cite{[R02]} has produced liftings in
the odd and even case to Witt vectors (while for the lifts we produce in
Theorem \ref{existminimallift}, their rationality cannot be
controlled). But his lifts are rarely minimally ramified.
His methods are those of Galois cohomology, while our lift is produced
using results of \cite{[T03]} and the formalism of deformation rings,
and thus our methods are quite indirect. Our methods work only in the
odd case.

In Part I of \cite{BoKh},  lifting  methods had also been used
to prove an analog
of Serre's conjecture for function fields.

The arguments of this paper use crucially the breakthroughs in
Wiles \cite{Wiles},
Taylor-Wiles \cite{TW}, and subsequent developments in modularity lifting
due to Taylor, Fujiwara, Diamond, Skinner-Wiles et al. The main idea of the paper leads to an
inductive approach to Serre's conjecture. There are 2 types of induction
involved, one on the number of primes ramified in the residual
representation
(see
Theorem \ref{levelone}), and the other on the residual characteristic
$p$ of the representation (see Theorem \ref{bertrand}). For the induction we
need a starting point and that is provided by results of Serre and Tate
which
prove the conjecture for level 1 representations mod 2 and 3. Given such a
starting point the method should in principle work for totally real fields.

We also use crucially in many places ideas or themes that
we have learnt from Serre's work.
His conjectures in \cite{Serre3} have been a great source of inspiration for
people in the field,
and at a more technical level his specification of the weight in
\cite{Serre3},
his  results
on relation between changing weight and $p$-part of the level, see
Th\'eor\`eme 11 of  \cite{Serre1},
his proof of the level one case of his conjectures for $p=3$
\cite{Serre2} page 710
which provides us
the toe-hold ({\it ongle de prise}!), and the theme
of studying compatible systems which was foregrounded via the many beautiful
results he proved about them have been a critical influence.

\vspace{3mm}

\noindent{\bf Acknowledgements:}
The first author would like to warmly thank Gebhard Bo\"ckle and Ravi
Ramakrishna through whom
he has learnt much about deformations of Galois representations, in the
course
of
his collaborations with them. He would also like to thank the
second author for the
invitation to visit Strasbourg, and the Department of Mathematics at
Strasbourg
for
its hospitality during the time when some of the work of this paper was
done.
We both would like to thank G.~B\"ockle, L.~Berger, C.~Breuil,
H.~Carayol,  R.~Schoof and J-P.~Serre for helpful
conversations/correspondence
in the course of this work.

\section{Minimal $p$-adic lifts of odd irreducible 2-dimensional Galois
 representations $\rhobar$}
\label{lift}

For $F$ a field, $\Q \subset F \subset \aQ$,
we write $G_{F}$ for the Galois group of $\aQ / F$.
For $\ell$ a prime,
we mean by $D_{\ell}$ (resp. $I_{\ell}$) a decomposition
(resp. inertia) subgroup of $G_{\Q}$ at  $\ell$.

Let $p$ an odd prime.
Fix  $\rhob: G_{\Q}\rightarrow \mathrm{GL}_2 (\aFp )$
to be an odd irreducible representation.
We assume that the Serre
weight $k(\rhob)$ is such that $2 \leq k(\rhob) \leq p+1$. (Note that there
is always
a twist of $\rhob$ by some power of the mod $p$ cyclotomic character
$\overline{\chi _p}$
that has weights in this range.)

Let $\F \subset \overline{\mathbb{F}_p}$ be a finite field
such that the image of $\overline{\rho}$ is contained in
$\mathrm{GL}_2 (\F )$, and let $W$ be the Witt vectors $W(\F )$.
By a {\it lift} of $\overline{\rho}$, we mean a continuous representation
$\rho : G_{\Q}\rightarrow
\mathrm{GL}_2 ({\cal O})$, where $\cal O$ is the ring of integers of a finite
extension of the field of fractions of  $W$, such that
the reduction of $\rho$ modulo the maximal ideal of ${\cal O}$ is
isomorphic to $\overline{\rho}$.

Let $\rho$ be such a lift and
let $\ell$ be a prime. One says that $\rho$ is {\it minimally
ramified at $\ell$} if it satisfies the following conditions:

- When  $\ell\not=p$,  it  is minimally
ramified at  $\ell$ in the terminology of \cite{[D97]}.
In particular, if $ \overline{\rho}$ is unramified at $\ell$,
$\rho$ is unramified at $\ell$. More
generally, when the image of $I_{\ell}$ is of order
prime to $p$, $\rho(I_{\ell})$ is isomorphic to its reduction
$\rhob (I_{\ell})$.

- When $\ell =p$ : If $k(\overline{\rho})\not= p+1$,
$\rho$ is minimally ramified at $p$ if $\rho $ is
crystalline of weights $(0,k (\overline{\rho} )-1 )$.
If $k(\overline{\rho})=p+1$,  $\rho$ is
{\it minimally ramified of   semi-stable type}  if $\rho $
is semi-stable non-crystalline of Hodge-Tate weights $(0, 1)$;
$\rho$  is {\it minimally ramified of crystalline
type} if $\rho $ is crystalline of Hodge-Tate weights
$(0, p)$.

\vspace{.5cm}

Let us make a few comments on the condition for $\ell =p$,
$k(\rhob )= p +1$.
Let $\chi _p: G_{\Q}\rightarrow \Z_p ^* $
be the $p$-adic cyclotomic character and
$\overline{\chi_p}$ its reduction modulo $p$.
If $k(\rhob)= p +1$, the restriction of $\rhob$
the  decomposition group $D_{p}$ is of the form:
$$\left( \begin{array}{cc}
{\overline{\chi_p}}\epsilon & \eta  \\
    0 & \epsilon
\end{array} \right),$$
where $\epsilon$ is an unramified character, and $\eta$
is a ``tr\`es ramifi\'e'' 1-cocycle, which corresponds
via Kummer theory to an element of $\Qp ^{*}\otimes \F$
whose image by the map defined by the valuation of $\Qp$
is a non-zero element of $\F$.

The lifting $\rho$ is minimally ramified of semi-stable type if
the restriction of $\rho$
to $I_{p}$ is of the form:
$$\left( \begin{array}{cc}
    \chi_p & *  \\
    0 & 1
\end{array} \right).$$
As Kummer theory
easily shows, this implies that the restriction of $\rho$ to the
decomposition
group $D_{p}$ is of the form:
$$\left( \begin{array}{cc}
    \chi_p \widehat{\epsilon} & *  \\
    0 & \widehat{\epsilon}
\end{array} \right),$$
where $\widehat{\epsilon}$ is an unramified character lifiting $\epsilon$.

The lifting $\rho$  is minimally ramified of crystalline type if
the restriction of  $\rho$
to $I_{p}$ is of the form:
$$\left( \begin{array}{cc}
    \chi_p ^{p}& *  \\
    0 & 1
\end{array} \right).$$

Indeed, by Bloch and Kato (3.9 of \cite{[BK]}), we know that such
representations are exactly the crystalline
non-irreducible representations of $D_p$
of Hodge-Tate weights $(0,p)$ . Furthermore, by
Berger, Li and Zhu (cor. 4.1.3. and prop. 4.1.4. of
\cite{BLZ}), we know that
the reduction of an irreducible crystalline representation
of $D_{p}$ of Hodge-Tate weights $(0,p)$ is isomorphic to
an unramified twist of
$\mathrm{ind}_{\Q _{p^2}}^{\Qp }(\omega _2 ^p)$, where $\Q _{p^2}$
is the quadratic unramified extension of $\Qp$ and $\omega_2$
is the fundamental character of level $2$: in particular, it is not
isomorphic
to the ``tr\`es ramifi\'ee'' representation.

The determinant of ${\rhob}$ is
$\overline{\chi_p}^{k(\rhob )-1}
\epsilon$ where $\epsilon$ is a character of conductor prime to $p$
(\cite{Serre3}). For $\ell\not=p$, the restriction to $I_{\ell}$
of the determinant of a minimal lift of $\rhob$ is the Teichmuller
lift (\cite{[D97]}). A semi-stable representation of $I_p$
of Hodge-Tate weights $(0,k-1)$ has determinant $\chi_p ^{k-1}$.
So we see that a minimal lift of $\rhob$ (of crystalline type
if $k(\rhob )=p+1$) has determinant
$\chi_p ^{k(\rhob )-1} \widehat{\epsilon}$, where
$\widehat{\epsilon}$ is the Teichmuller lift of $\epsilon$.
If $k(\rhob )=p+1$, a minimal lift of semi-stable type has determinant
$\chi_p  \widehat{\epsilon}$.

As it is suggested in 5.2. of \cite{[KR03]},
using B\"ockle's appendix to \cite{[K03]} and
Taylor's theorems in \cite{[T03]} and \cite{[T02]},
one can produce minimal lifts of Serre type $\rhob$:

\begin{theorem}\label{existminimallift}
Let $p$ be a prime $> 3$.
Let
$\overline{\rho} : G_{\Q}\rightarrow
\mathrm{GL}_2 (\F)$ be an odd absolutely irreducible representation.
We suppose that $2 \leq k(\rhob) \leq p+1$ and 
 $k(\rhob)\not= p$. Then
$\rhob$ has  a lift $\rho$ which is minimally ramified
at every $\ell$, and if the Serre weight is
$k(\overline{\rho} )=p+1$,  one can impose
that $\rho$ be either of crystalline type (of weight $p+1$)
or of semi-stable type (of weight $2$).
\end{theorem}

The rest of the section will be devoted to the proof of Theorem
\ref{existminimallift}. We deal with the cases in Theorem
\ref{existminimallift} when  the image of $\rhob$ is
solvable right away. In this case the asserted lifts in Theorem
\ref{existminimallift} come from the fact that by Langlands-Tunnell one
knows the modularity of $\rhobar$ and then one uses the fact that the
recipe in \cite{Serre3}
for optimal weights and levels from which $\rhob$ arises  has been proven to
be correct
as a result of the work of a number of people, see \cite{Ribet1}
(note that we are
assuming that $p>3$). In the weight $p+1$ case we see that there is a
semistable weight 2 lift (that arises from a newform in fact)
which is minimal,
starting from the crystalline weight $p+1$ lift (which we know arises
from
a newform),
by using a standard argument
that relies on the fact that $\F$-valued functions on the  projective
line
over $\F_p$ as a $GL_2(\F_p)$ module decomposes as a sum of the
trivial representation and ${\rm Symm}^{p-1}(\F_p^2)$.

From now on, we suppose that the image of $\rhob$ is not solvable.
For the proof of the theorem, we have to consider minimally
ramified deformations
$G_{\Q }\rightarrow \GL_ 2 (R)$ of $\overline{\rho }$.
Let $R$ be a local profinite $W$-algebra, with an isomorphism
of $R/\mathcal{M}_R$ with $\F$ with $\mathcal{M}_R$ the maximal ideal
of $R$ ($W$ is as above the Witt ring $W(\F )$).
A deformation of $\rhob$ is a continuous representation
$\gamma :G_{\Q}\rightarrow \GL_2 (R)$
such that $\gamma$ $\mathrm{mod} \mathcal{M}_R$
is ${\rhob}$, where we take $\gamma$ up to conjugation by matrices that
are 1 mod $\mathcal{M}_R$. We say that the deformation is
minimally ramified, if:

- for $\ell \not=p$, $\gamma$ is minimal in the sense of \cite{[D97]};

- if $k(\rhob )<p$, the restriction of $\gamma$ to $D_p$ comes
from a Fontaine-Laffaille module;

- if  $k(\rhob )=p+1$, the restriction of $\gamma$ to $I_p$
is of the form:
$$\left( \begin{array}{cc}
    \chi_p ^{k-1} & *  \\
    0 & 1
\end{array} \right),$$
with $k=p+1$ if we are in the crystalline type, and $k=2$ if we are
in the semi-stable type.

The condition of being minimally ramified is a deformation condition in
the sense of \cite{[M97]}, and hence the minimally ramified deformation
problem
has a universal object. More precisely,
if $k(\rhob) \not= p+1$, there exists a universal
minimally ramified deformation
$\rho_{\mathrm{univ}}: G_{\Q }\rightarrow \GL_2 (R_{\mathrm{univ}})$;
if $k(\rhob )=p+1$, we have two universal rings
$R_{\mathrm{univ,ss}}$ and $R_{\mathrm{univ,crys}}$.

Theorem \ref{existminimallift} follows from:
\begin{theorem}\label{completeinters} Let ${\rhob}$ as in Theorem
\ref{existminimallift}. Suppose that the image of
$G_{\Q }$ is not solvable. Then,  $R_{\mathrm{univ}}$ if
$k({\rhob} )\not= p+1$, and
$R_{\mathrm{univ,crys}}$ and $R_{\mathrm{univ,ss}}$ if
$k({\rhob} )= p+1$, are finite
flat complete intersection $W$-algebras. \end{theorem}

\begin{proof}
Define for each $\ell$, the $W$-algebra  $R_{\ell }$ of versal
deformations of $\overline{\rho}_{\mid D_{\ell }}$ which are
minimally ramified (if $k({\rhob} )=p+1$,
we have to consider the two $W$-algebras
$R_{p,\mathrm{crys} }$ and $R_{p,\mathrm{ss} }$)
and such that the determinant is the restriction
to $D_{\ell}$ of $\chi_p ^{k-1} \widehat{\epsilon}$,
with $k=k(\rhob )$ except in the case $k(\rhob )=p+1$
and we are in the case of semi-stable type, and then $k=2$.

Proposition 1 of B\"{o}ckle (appendix to
the article of the first author \cite{[K03]}) says that
if
the $W$-algebras  $R_{\ell }$ are flat,  complete intersections of
relative dimension $\mathrm{dim}_{\kappa }(H^0 (D_{\ell },
\mathrm{ad}^0 (\overline{\rho }))+\epsilon_{\ell }$ with
$\epsilon_{\ell}=0$ if $\ell \not= p$ and $\epsilon _{p }=1$, then
the $W$-algebra $R_{\mathrm{univ}}$ has a presentation as a CNL
$W$-algebra as $$W[[X_1,\cdots, X_r]]/(f_1,\cdots,f_s)$$ with $r \geq s$.
Recall that  $\mathrm{ad}^0 (\overline{\rho })$ is the subspace
in the adjoint representation of matrices of trace $0$.

Except in the case of $k({\rhob} )=p+1$
and $R_{\ell}$ is $R_{p,\mathrm{crys} }$, it is proved by Ramakrishna
(\cite{[R02]}) and Taylor (\cite{[T03ico]}) that
$R_{\ell}$ is smooth over $W$ of the  dimension
asked for in B\"{o}ckle's proposition.
For $R_{\ell}=R_{p,\mathrm{crys} }$ in the case when $k({\rhob}
)=p+1$,
it is proved by
B\"{o}ckle that $R_{p,\mathrm{crys} }$ is  a relative complete
intersection of relative dimension $1$ (remark 7.5 (iii)
of \cite{[B00]}). For the convenience of the reader we give a proof of
this result of B\"ockle.

\begin{prop}\label{p+1} In the case $k(\rhob)=p+1$, the $W$-algebra
$R_{p,\mathrm{crys} }$ is formally
smooth of dimension $1$. \end{prop}

\begin{proof}
Let $U$ be the the $\F$-vector space underlying the representation
$\rhob$. Call  $F^1 U$ the line stable under $D_p$.
Let $F^0 $ be the sub-space of the endomorphisms of
$ \mathrm{ad}^0 (\rhob )$  which respects the filtration $U\supset F^1 U$
and
let $F^1 \subset F^0$ be the sub-space of elements of $F^0$
that act trivially on $U/F^1 U\oplus F^1 U$. As deformations
of crystalline or semi-stable type are easily seen to be triangular,
we see that the relative cotangent spaces $\mathcal{M}/(p,\mathcal{M}^2)$
for the crystalline or semi-stable type deformations are isomorphic to
the kernel of the map :
$$H^1( D_p, F^0 )\rightarrow
H^1( I_p, F^0 / F^1 ) ^{D_p / I_p}.$$
One knows that, as $\rhob$ is ``tr\`es ramifi\'ee'', that the dimension
of the kernel is $1$ (Lemma 29 of \cite{[DDT]}).
We see that $R_{p,\mathrm{crys} }$ is a quotient of $W[[T]]$.
If $f$ is a non zero element of $W[[T]]$, there are only finitely many
morphisms of $W$-algebras $W[[T]]/(f) \rightarrow W$.
Thus to prove that $R_{p,\mathrm{crys} }$ is isomorphic to $W[[T]]$,
it suffices to show  that
$\rhob$ has infinitely many  inequivalent liftings $D_p \rightarrow
\mathrm{GL}_2 (W )$.
Let $n\geq 1$ be an integer.
Let $\rho _n : D_p \rightarrow \mathrm{GL}_2 (W/ p^n )$
be a lifting of $\rhob$ of crystalline type.
Let us lift it modulo $p^{n+1}$.

In a convenient basis,
$\rho _n$ is of the form :
$$\left( \begin{array}{cc}
    \delta & \eta  \\
    0 & 1
\end{array} \right),$$
where  $\delta : D_p \rightarrow (W/p^n )^* $ is a character
whose reduction modulo $p$ is $\overline{\chi _p}$ and
whose restriction to $I_p$ is the reduction mod. $ p^n$
of $\chi_p ^p$.
Let $\widehat{\delta}$ be a lifting mod. $p^{n+1}$ of
$\delta$ whose restriction to $I_p$ is the reduction
of  $\chi_p ^p$. Let $\gamma : D_p \rightarrow \F$
be an unramified character.  We define the character
$\delta _{\gamma} : D_p \rightarrow
(W/p^{n+1} W)^*$ by :
$\delta _{\gamma}(\sigma )=1+p^n \gamma (\sigma )\
\mathrm{mod} p^{n+1}$.

Write $f_{\gamma }$ for the connecting homomorphism :
$ H^1 (W/ p^n W (\delta ))\rightarrow H^2 (\F (\overline{\chi _p } ))$
for the exact sequence :
$$(0) \rightarrow \F (\overline{\chi _p } ) \rightarrow
W/ p^{n+1} W ( \widehat{\delta}
\delta _{\gamma})\rightarrow W/p^{n}W ( \delta)\rightarrow (0)$$
(the cohomology is for the group $D_p$).
As the restriction to $I_p$ of the reduction
modulo $p^2$ of $\chi _p \ \widehat{\delta}^{-1}
\delta _{\gamma}^{-1}$ is non trivial, the map :
$$ H^0 (W/ p^{n+1} W (\chi _p \ \widehat{\delta}^{-1}
\delta _{\gamma}^{-1}))\rightarrow H^{0} (\F )$$
is zero. By Tate local duality this implies that $f_{\gamma }$ is surjective. We have
the exact sequence :
$$ (0)\rightarrow H^1 (\F (\overline{\chi _p } ))\rightarrow
H^1 (W/ p^{n+1} W ( \widehat{\delta}
\delta _{\gamma}))
\rightarrow
H^1 (W/ p^{n} W ( \delta))\rightarrow H^2 (\F (\overline{\chi _p } ))
\rightarrow (0).$$

Let us still denote by $\eta$ the  class
in $H^1(\F (\overline{\chi _p } ) )$ of the cocycle $\eta$.
A direct calculation gives that
$f_{\gamma}(\eta)- f_0 (\eta )$ is the cup product of
$\gamma\in H^1 (\F )$ and the reduction
$\overline{\eta }\in H^1 (\F (\overline{\chi _p}))$
of $\eta$.
As $\overline{\eta }$ is "tr\`es ramifi\'e", class field theory
implies that, if $\gamma$ is non zero, this cup product is non zero.
We find a (unique) $\gamma_1$ such that the
$f_{\gamma_1}(\eta)=0$. This implies that
$\eta$ lifts to an $\widehat{\eta}\in H^1 (W/ p^{n+1} W ( \widehat{\delta}
\delta _{\gamma_1}))$. The representation
$\rho _n$  has
a crystalline type lifting modulo $p^{n+1}$.
The cardinal $\mid  H^1 (\F (\overline{\chi _p } )) \mid$
is $\mid \F \mid ^2$. As two cohomology classes $\widehat{\eta}$ give
rise to equivalent liftings exactly when they differ by an element
of $1+p^n W/1+p^{n+1} W$,  we have $\mid \F\mid$
non equivalent liftings of crystalline type modulo $p^{n+1}$ of
$\rho _n$. This proves that $\rhob$ has infinitely many liftings
to $W$ and proves the proposition.
\end{proof}

From now on we denote by $R_{\mathrm{univ}}$ the deformation
ring we consider.
If we can prove that $R_{\mathrm{univ}}$ is a finitely generated
$W$-module, we are done by a standard argument (see for example Lemma 2
of the above quoted appendix of B\"{o}ckle). Namely we see easily
that the sequence $f_1,\cdots,f_s,p$ is regular, and this gives
that $R_{\mathrm{univ}}$ is a finite flat complete intersection over $W$.
So one has only to prove that $R_{\mathrm{univ}}/pR_{\mathrm{univ}}$
is of finite cardinality. We prove a lemma that reduces this to proving that
$\rho_{\mathrm{univ}}$ mod $p$,
which we denote by $\overline{\rho_{\mathrm{univ}}}$, has finite image
(this is inspired
by Lemma 3.15 of de Jong (\cite{[dJ01]}):

\begin{lemma}\label{finite} Let $\kappa$ be a finite field of characteristic
$p$, $G$ a profinite group satisfying the $p$-finiteness
condition (chap. 1 of Mazur \cite{[M97]})
and $\eta: G\rightarrow \GL _N (\kappa )$ be an
absolutely irreducible continuous representation. Let
$\mathcal{F}_N (\kappa )$ a subcategory of  deformations of $\eta $
in $\kappa$-algebras which satisfy the conditions
of 23 of \cite{[M97]}. Let
$\eta_{\mathcal{F}} : G\rightarrow \GL_N (R_{\mathcal{F}} )$
be the universal deformation
of $\eta$ in $\mathcal{F}_N (\kappa )$.
Then $R_{\mathcal{F}}$  is finite if and only
if $\eta_{\mathcal{F}} (G )$ is
finite.\end{lemma}

\begin{proof} It is clear that if $R_{\mathcal{F}}$
is finite,  $\eta_{\mathcal{F}} (G )$ is finite.
Let us suppose that $\eta _{\mathcal{F}}(G )$ is finite. As
$\eta$ is absolutely irreducible, a theorem of Carayol
says that $R_{\mathcal{F}}$ is generated by the traces of
the $\eta_{\mathcal{F}}(g)$, $g\in G$ (\cite{[C94]}).
As $\eta _{\mathcal{F}}(G )$ is finite,
for each prime ideal $\wp$ of $R_{\mathcal{F}}$,
the images of these traces in
the quotient $R_{\mathcal{F}}/ \wp$
are sums of roots of unity, and  there is a finite number of them.
We see that  $R_{\mathcal{F}}/ \wp$ is a finite extension of
$\kappa$. It follows that the noetherian ring
$R_{\mathcal{F}}$ is of dimension $0$, and so is finite.\end{proof}

We show now that
$\overline{\rho_{\mathrm{univ}}}$
has finite image.
We begin by proving (see also Lemma 2.12 of \cite{[dJ01]}):

{\it Claim:} for each $\ell\not= p$,
$\overline{\rho_{\mathrm{univ}}}$ is
finitely ramified at $\ell$. In fact, the order of
$\overline{\rho_{\mathrm{univ}}}(I_\ell)$ is the same as that of
$\rhobar(I_\ell)$.

{\it Proof of claim:}
The only case that needs argument is when
the restriction of $\rhobar$
to $I_{\ell}$ is of type:
$$ \xi \otimes \left( \begin{array}{cc}
     1 & \phi  \\
     0 & 1
 \end{array} \right),$$
with $\phi$ a ramified character. The minimality condition
implies that the restriction of  ${\rho}_{\mathrm{univ}}$
to $I_{\ell}$ is of the form:
$$ \tilde{\xi} \otimes \left( \begin{array}{cc}
     1 & \tilde{\phi}  \\
     0 & 1
 \end{array} \right),$$
with $\tilde{\xi}$ beeing the Teichmuller lift of $\xi$.
The morphism $\tilde{\phi}$ is tamely ramified, so its image is cyclic.
As $R_{\mathrm{univ}}/pR_{\mathrm{univ}}$ is a $\F_p$-algebra,  $p
\tilde{\phi}=0$ and
$\tilde{\phi}$ has image of order $p$. This proves the claim.

\vspace{3mm}

We have the following crucial proposition that follows easily from
the results of Taylor in \cite{[T02]} (in the ordinary case) and \cite{[T03]} (in
the supersingular case), and
is the key input in the proof of Theorem \ref{completeinters}.
For a totally real field $F$ of even degree, let $B_F$ be the unique  (up to
isomorphism) totally definite quaternion algebra that is unramified at
all finite places. We recall that we are assuming that the image of
$\rhob$ is not solvable.

\begin{prop}\label{taylor}
There exists a totally real field $F$ of even degree with the following
properties:

(i) $F$ is unramified at $p$, and in the case when $\rhob|_{D_p}$
is irreducible, $F$ is split at $p$,

(ii) $\rhobar|_{G_F}$ is absolutely irreducible and non-dihedral,

(iii) $\rhob|_{G_F}$ is unramified
outside the primes above $p$ and arises from
a cuspidal automorphic form $\pi$ for $B_F$ that is unramified at all
finite places and is of weight $k=k(\rhobar)$ at all the infinite
places, and which is ordinary at places above $p$ when $\rhob$ is ordinary
at $p$
($\rhob|_{I_p}$ has a quotient of dimension $1$ with trivial
action).
In the weight $p+1$ case there is also a cuspidal automorphic form
$\pi'$ for $B_F$ that is unramified at all
finite places prime to $p$, is the Steinberg representation at all
places above $p$, and is of weight $2$ at all the infinite
places, and which is ordinary at places above $p$.
\end{prop}

\begin{proof}
Property (ii) is satisfied for all totally real $F$ by lemma
\ref{potirr}.

The supersingular case is covered in
\cite{[T03]}
explicitly (see Theorem 5.7 of
\cite{[T03]}), while the ordinary case is not but can be deduced from the
arguments of \cite{[T02]}.
Thus we treat below only the case when $\rhob$ is ordinary.

Let us borrow for a moment the
notations of Taylor \cite{[T02]}, even if it contradicts the notations
of this paper,
for the next 4 paragraphs.

Let $A$ be the  abelian variety given by application
of Moret-Bailly's theorem in \cite{[T02]} p. 137
(so we make the twist at the bottom of p. 136).
In particular, $A$ is defined over the
real field $E$ (our $F$), $\rhob$ has values in a finite field
of characteristic $\ell$ (our $p$). The
abelian variety $A$ is of
Hilbert-Blumenthal type with multiplication by the real field $M$.
There is a prime $\lambda$ of $M$ above $\ell$ such that the
restriction of $\rhob$ to $G_E$ is  isomorphic to the $G_F$-
representation on  $A[\lambda]$, the points of $A$ killed
by $\lambda$. The compatible system of $G_E$-representations
attached to $A$ is modular, say comes from an
automorphic form
$\pi _A$ of parallel weight $2$.

Let $x$ a place of $E$ above $\ell$. We use Lemma 1.5
of \cite{[T02]} to get the needed information for  $(\pi _A )_x$,
namely we will prove that it is ordinary (with respect to
the place $\lambda$ of $M$).
Let $n=\ell-k(\rhob )+1$ if $k(\rhob)\not=2$ and
$n=0$ if $k(\rhob)=2$. Note that $n$ is as in Lemma 1.5 of
\cite{[T02]}.
Note that as we are assuming $k(\rhob)\not= \ell$,
we have $n\not=1$ and Lemma 1.5  applies.
We have $0\leq n< \ell-1$ and for a place $x$ of
$E$ above $l$, the $\lambda$-adic representation arising from $A$
when restricted  of the decomposition group
$G_x$ is of the form:
$$\left( \begin{array}{cc}
    \epsilon \chi _1 & *  \\
    0 & \chi_2
\end{array} \right),$$
with $\chi_2$ unramified and the restriction of $\chi_1$
to the inertia subgroup $I_x$ of $G_x$ is the reduction of
$\epsilon^{-n}$, $\epsilon$ being the cyclotomic character
(acting on $\ell^*$ roots of unity). We know by the proof
of Lemma 1.5 that $A$ has mutiplicative reduction over $E_x$
or good reduction over $E_x (\zeta_{\ell })$. Furthermore, there
is a prime $\wp$ above $p\not=\ell$ such that the action
of $G_x$ on $A[\wp ]$ has the form $\psi _1 \oplus \psi _2$,
with $\psi_2$ unramified and the restriction of $\psi _1$
to $I_x$ is $\omega^{-n}$ where $\omega$ still denote
the reduction of the Teichmuller lift $\omega$ of $\epsilon$.

In the case $n=0$ ($k(\rhob )=\ell +1$ or 2),
we see by looking at the Tate module $T_{\wp} (A)$
that $A$ has semistable ordinary reduction over $E_x$.
If $k(\rhob )=\ell +1$, $A$ has mutiplicative reduction
at all $x$ over $\ell$
and $(\pi_{A})_{x}$ is Steinberg.
When  $k(\rhob )=2$, and $\chi_{1}\chi_{2}^{-1}\not= 1$,
Taylor finds, for a place $v$  of $F$ above $\ell$,
an abelian variety $A_{v}$ over $F_{v}$ with ordinary
good reduction.
The theorem of Moret-Bailly  \cite{Moret}
produces for us
an abelian variety $A$ with good reduction at all primes $x$
of $E$ above $\ell$ such that the
restriction of $\rhob$ to $G_E$ is  isomorphic to the $G_F$-
representation on  $A[\lambda]$. We see that, if we choose
$A$ like this, $(\pi_{A})_{x}$ is unramified.

(If $k(\rhob )=2$ and $\chi_{1}\chi_{2}^{-1}= 1$, we are in the
case $\chi_{v}^2 =1$ of the proof of Lemma 1.2 of Taylor. But, as the
restriction of $\rhob$ to $G_{v}$ comes from
a finite flat group scheme over the ring of integers of $F_{v}$,
we can choose the abelian variety
$A_{v}$ that figures in Lemma 1.2 which good ordinary reduction,
by the same arguments as Taylor uses when
$\chi_{v}^2 \not=1$. This is because the class in
$H^1 (G_{v}, O_{M}/\lambda (\epsilon ))$
of the extension defined by $\rhob _{\mid G_{v}}$ comes
from units by Kummer theory. Then, as above, we can
choose $A$ with good ordinary reduction at all places $x$ of
$E$ above $\ell$. Then $(\pi _{A})_{x}$ is unramified at these places.
We have proved the proposition if $k(\rhob )=2$.
If we proceed this way, we  may obviate the use of level-lowering results of
\cite{Jarvis} at the end of the proof of the proposition.)

Suppose now $n\not=0$. Then, looking at the
Tate module $T_{\wp}(A)$, we see that the abelian variety $A$
has good reduction
over $E_x (\zeta_{\ell })$. Let $A[\lambda]^0$ and
$A[\lambda]^{\mathrm{et}}$ the connected and etale components
of the $\lambda$-kernel of the reduction at $x$ of the N\'eron model
of $A$ over $E_x (\zeta_{\ell })$.
Let $T_{\lambda } (A)$,
$T_{\lambda }^0 (A)$
and $T_{\lambda }^{\mathrm{et}} (A)$ be the corresponding Tate-modules and
let $D$, $D^0$ and $D^{\mathrm{et}}$ the correponding Dieudonn\'e modules.
We have $D=D^0 \oplus D^{\mathrm{et}}$.
Taylor proves in  Lemma
1.5 that $I_x$ acts on
$\mathrm{Lie}(A[\lambda]^0 )$ by multiplication by $\omega^{-n}$ and
trivially
on $A[\lambda]^{\mathrm{et}}$.
As the action of $I_x$ on $D$ factors through $\mathrm{Gal}(E_x (\zeta_{\ell
})
/E_x )$,
it follows that $I_x$ acts on $D$ by multiplication
by $\omega^{-n}$ on $D^0$ and trivially on $D^{\mathrm{et}}$.
By the appendix B in Conrad-Diamond-Taylor \cite{CDT} it follows that the
action
of the Weil-Deligne group $WD_x$ on the compatible system
of Galois representations attached to $A$ factors through
the Weil group and has the form $\eta_1\oplus \eta_2$, with
$\eta_2$ unramified and $\eta_2 (\mathrm{Frob}_x )$ a $\lambda$-adic
unit, and the restriction of $\eta_1$ to $I_x$ is $\omega ^{-n}$.
It follows that
$(\pi_A )$ is ordinary at $x$, with nebentypus $\Psi$
such that $\Psi \omega^{-n}$ is unramified
at every place of $E$. Using the base change technique of Skinner-Wiles
\cite{SW1} we may also ensure, after base change to
 a totally real solvable extension of $E$ that is unramified at primes
 above $p$, that $\pi_A$ is unramified at primes
not lying above $p$ (and is still ordinary at primes above $p$).

We revert now to the notation of the present paper, i.e., $\ell$ is now
$p$ and $E$ is $F$.

We use weight shifting arguments due to Hida which
give that projection onto the highest weight vector of coefficients
induces an isomorphism on the ordinary part of the cohomology. Namely, by
the arguments in Section 8 of \cite{Hida1} (note that his arguments apply in our
situation,
as $p>3$ and $p$ is unramified in $F$ and then
use Lemma 1.1 of \cite{[T03]}, although the neatness assumption of
\cite{Hida}
need not be satified here)
we may deduce from the previous paragraph (i.e., the existence of
ordinary $\pi_A$ as above that is of weight $(2,\cdots,2)$ at
infinity, unramified away from $p$, and at places $\wp$ above $p$ is
principal series of conductor dividing $\wp$)
that $\rhob$ arises from a cuspidal automorphic representation $\pi$ of
$B_F(\A_F)$
which is
is unramified at all
finite places that are prime to $p$, at places $\wp$ above $p$ is
either an unramified
principal series or Steinberg (and hence again of conductor dividing $\wp$),
and is of weight $(k(\rhob),\cdots,k(\rhob))$ at infinity.
Further $\pi$ is ordinary
for all such $\wp$.  Also
when $k(\rhob)=p+1$, one may choose a $\pi'$ as above except that it
has weight $(2,\cdots,2)$ at infinity.
Now when $k(\rhob)$ is not $2$ we are done as
then forms that are Steinberg at a place above $p$ cannot be ordinary
at that prime for weights bigger than 2. The
Serre weight 2 case needs an additional argument.
We have to ensure that one can
choose $\pi$ of parallel weight $(2,\cdots,2)$, which is unramified
at all places and gives rise to $\rhob$. This follows from arguments
using Mazur's principle proved by Jarvis (see Theorem 6.2 of \cite{Jarvis}).

\end{proof}

We prove the general well-known lemma used in the proof above.

\begin{lemma}\label{potirr} ($p>2$)
Let $\eta : \Gal \rightarrow \mathrm{GL}_2 (\aFp )$ be an odd
Galois representation. Suppose that the image of $\eta$
is not solvable (so $\eta$ is absolutely irreducible). 
Let $F'$ be a Galois totally real
finite extension of $\Q$ contained in $\aQ$ and $F''$ be a quadratic
extension of $F'$. Then  $\eta (G_{F''})$
is not solvable ; in particular, the restriction 
of $\eta$ to $G_{F'}$ is non dihedral.
\end{lemma}
\begin{proof}
Let $H$ be the the image of $G_{\Q}$ in $\mathrm{PGL}_2 (\aFp )$.
So, by Dickson (th. 2.47. of \cite{[DDT]}), $H$ is conjugate to
$\mathrm{PSL}_2 (\mathbb{F}_{p^r})$ or
$\mathrm{PGL}_2 (\mathbb{F}_{p^r})$ for $p^r\not= 2,3$,
or is isomorphic to $A_5$ : 
the $p^r = 2,3$, triangular, dihedral, $A_4$, $S_4$,
cases are excluded as  $H$ is not solvable.

As $\eta$ is odd, the image of a complex conjugation in $H$ is non trivial.
The image $H'$ of $G_{F'}$ in
$H$ is a non trivial normal subgroup. 
If $H$ is $\mathrm{PSL}_2 (\mathbb{F}_{p^r})$ or
$\mathrm{PGL}_2 (\mathbb{F}_{p^r})$ with $p^r\not= 2,3$, $H'$ contains
the simple group $\mathrm{PSL}_2 (\mathbb{F}_{p^r})$.
In the $A_5$ case, $H'=A_5$. 
As $\mathrm{PSL}_2 (\mathbb{F}_{p^r})$ and $A_5$ do not have 
a subgroup of indice 2, the image of $G_{F''}$ in 
$\mathrm{PGL}_2 (\aFp )$ contains $\mathrm{PSL}_2 (\mathbb{F}_{p^r})$
or $A_5$. This proves
the lemma.\end{proof}

We return to the proof of Theorem \ref{completeinters}.
Choose $F$ as in Lemma \ref{taylor}.
We show that $\overline{\rho_{\mathrm{univ}}}|_{G_F}$
has finite image for this choice of $F$, which clearly implies that
$\overline{\rho_{\mathrm{univ}}}$
has finite image, which by Lemma \ref{finite}
implies Theorem \ref{completeinters}.

Let $\rho_{\mathrm{univ},F}:G_{F}\rightarrow \GL_2
(R_{\mathrm{univ},F})$ be the universal,
minimally ramified
$W$-deformation of the restriction of $\overline{\rho}$ to $G_F$:
recall that it is unramified at every prime of $F$ of residual
characteristic
$\not= p$ and is minimally ramified at every prime above $p$, and for
primes above $p$, we take the same conditions
as we have taken to define $\rho_{\mathrm{univ}}$ (and their variants
for the 2 deformation rings when $k(\rhobar)=p+1$).
Because of the claim we proved above
$\overline{\rho_{\mathrm{univ}}}|_{G_F}$
is a specialisation of $\overline{\rho_{\mathrm{univ},F}}$, amd it
will be enough to prove that
$\overline{\rho_{\mathrm{univ},F}}$ has finite image.

From Fujiwara's generalisation in \cite{Fujiwara}
of $R=T$ theorems of \cite{Wiles} to the case of totally real fields
 in the ordinary case (see also \cite{SW3}),
and Theorem 3.3 of
\cite{[T03]} in the supersingular case (note that we can apply these lifting
theorems because of Lemma \ref{potirr}, and as we are
excluding weight $p$, and in weight $p+1$ all our lifts are ordinary)
we deduce from Proposition \ref{taylor}
that $R_{\mathrm{univ},F}$ is finite as a $W$-module. 
This comes from
proving an $R=T$ theorem for $\rhob|_{G_{F}}$,
which identify the universal deformation ring with the
completion of the ordinary Hecke algebra of Hida acting on cusp forms over
$F$
of level
1 and parallel weight $k(\rhob)$ (except when we consider
the ``tr\`es ramifi\'e'' case and the  deformation  ring
$R_{\mathrm{univ,ss}}$ when the
cusp forms will be of weight 2 unramified at all primes not above $p$,
and
for primes above $p$ will be special). From this it follows that
$\overline{\rho_{\mathrm{univ},F}}$ has finite image, and hence we are
done
with the proof of Theorem \ref{completeinters}.

\end{proof}

\section{Minimal compatible systems that lift semistable $\rhobar$}

Let $p$ be  prime $>3$, and $\rhob$ be as before of Serre type,
i.e., an
odd, irreducible 2 dimensional representation of $\Gal$. We further assume
from
now on that it is {\it semistable}.
By semistable we mean as usual that for primes outside
$p$ ramification is unipotent, and that at $p$ (slightly unusually)
we require that
the Serre weight $k(\rhob)$ is $\leq p+1$ (all
representations after
twisting by
some power of the mod $p$ cyclotomic character ${\overline{\chi_p}}$ have
such a
weight). We remark that the determinant of such a $\rhob$
is $\overline{\chi_{p}}^{k(\rhob )-1}$, hence $k(\rhob )$ is even.
As $k(\rhob )\not= p$, Theorem \ref{existminimallift} applies.


The following theorem is easily deduced from the arguments of
Dieulefait, see \cite{D}
and \cite{Wint}, after Theorem \ref{existminimallift}.
The proofs of \cite{D} and \cite{Wint}
relies on the method of proof of Theorem 6.6 of
Taylor's paper \cite{[T03]}.

Below, for a
2-dimensional  $p$-adic representation to be crystalline at $p$ of weight
$k$
 we require that the Hodge-Tate weights are
$(0,k-1)$. By strictly
compatible system of representations of $\Gal$,
we mean that for a given prime $q$ the
$F$-semisimplification of the Weil-Deligne
parameter at $q$ is isomorphic for all the places $\lambda$
(including places $\lambda$ whose residue characteristic $\ell$ is $q$), while by compatible we require this
only for
$\lambda$ whose residue characteristic is prime to $q$.

\begin{theorem}\label{minimal} Let $\rhobar$ be a semistable representation
of Serre type. Then $\rhobar$ has a
minimally ramified lift $\rho$ by Theorem \ref{existminimallift} which
at $p$ is either crystalline of weight $k(\rhob)$,
or in the weight $p+1$ case can also be chosen to be of
semistable type of weight 2 at $p$. We can further impose :

(i) assume $\rho$ is unramified outside $p$ and crystalline at $p$
(equivalently, $\rhobar$ is of level 1 and in the weight $p+1$ case we
consider the crystalline lift). There
is a number field $E$ and a
compatible system of representations $(\rho_\lambda)$ with
$\lambda$ running through the  set of finite places of $E$
such that at a place above $p$, the member of the
compatible system at that place is $\rho$. Further the Weil-Deligne
parameters are unramified at all primes except
perhaps for $\lambda$ and $q$ of the same characteristic 2 (thus, in
particular,
$\rho_{\lambda}$
is crystalline at places whose residue characteristic is the same as
that of $\lambda$ and is not 2).

(ii) Assume  the condition of (i) is not satisfied, thus either $\rho$ is
ramified at a prime outside $p$, or at $p$ it is semistable of weight
2 (equivalently, $\rhobar$ is not of level 1, or in the weight $p+1$ case we
consider the semistable weight 2 lifting).  Then there
is a number field $E$ and a strictly
 compatible system of representations $(\rho_\lambda)$ with
$\lambda$ running through the  set of finite places of $E$
such that at a place above $p$, the member of the
compatible system at that place is $\rho$. Further $(\rho_{\lambda})$
arises from the \'etale cohomology of some variety over $\Q$, and when
$k(\rhobar)=2$ or $p+1$ (and we consider weight 2 semistable lifting
$\rhobar$), we may choose $E$ so that there exists an abelian
variety $A$ over $\Q$ of dimension $[E:\Q]$ and an embedding ${\cal
 O}_E \hookrightarrow {\rm End}(A/\Q)$ such that $(\rho_{\lambda})$
is equivalent to the representation on the $\lambda$-adic Tate module
of $A$. Further $A$ has semistable and bad (multiplicative) reduction
reduction only at the primes dividing the prime to $p$ part of the Artin
conductor of $\rhob$, and also $p$ when $k(\rhob)=p+1$.
\end{theorem}

\begin{proof}

We sketch the proof for completeness, although as said above this is
proved in \cite{D} and \cite{Wint}.

When $\rhobar$ has
solvable image this follows from the fact that Serre's conjecture is
known for $\rhobar$.

So we assume that the image of $\rhob$ is not solvable.
Then, by Theorem 3.3 and lemma 5.6. of \cite{[T03]} (supersingular case) and
\cite{[T02]}, \cite{Fujiwara}, \cite{SW3} (ordinary case),  we know that
there is a totally real field $F$ unramified at $p$, of even degree,
such that $\rho|_{G_F}$
arises
from a cuspidal automorphic representation $\pi$ of $GL_2(\A_F)$
that is holomorphic of parallel weight
$(k(\rhobar),\cdots,k(\rhobar))$ (or also from such a $\pi$ of weight
$(2,\cdots,2)$ when $k(\rhob)=p+1$), and such that the local components of
$\pi$ at finite
places are either Steinberg or unramified principal series.

Furthermore, by Arthur-Clozel solvable base change (\cite{AC}),
we know that for each $F'\subset F$ such that $F/ F'$ has solvable
Galois group, the restriction of $\rho$ to $G_{F'}$ comes
from an automorphic representation $\pi _{F'}$ of $GL_2(\A_{F'})$.
By the argument using Brauer's theorem as in proof of Theorem 6.6 of
\cite{[T03]},
we then have a finite extension $E$ of $\Q$ and a system
$(\rho _{\lambda})$ of irreducible representations of
$\mathrm{Gal}(\aQ / \Q)$,
$\lambda$ describing the finite places of $E$.
There is a place $\lambda$ of $E$ above $p$ such that $\rho_{\lambda }$
is isomorphic to $\rho$.
The system
$(\rho _{\lambda})$  satifies
the weak compatibility property that there exists a finite
set $S$  of primes of $\Q$
such that $\rho_{\lambda}$ is unramified
outside $S$ and the residual characteristic $\ell$ of
$\lambda$, and such that, for $q\notin S$, the characteristic polynomials
of $\rho_{\lambda}(\mathrm{Frob}_{q})$ are the same for
$\lambda$ not over $q$.

To get the finer compatibility properties, one uses that, for
$F'\subset F$ such that $F/ F'$ has solvable Galois group,
the system $(\rho _{\lambda})$ restricted to $G_{F'}$ comes from $\pi_{F'}$. Let
$q$ be a prime number. Let $\mathcal{Q}$ be a prime of $F$
above $q$ and let $F(\mathcal{Q})$ be the subfield
of $F$ fixed by the decomposition group
$\subset \mathrm{Gal}(F/\Q )$ at $\mathcal{Q}$.
We know that the restriction of $(\rho_{\lambda})$
to $G_{ F(\mathcal{Q})}$ comes from $\pi_ { F(\mathcal{Q})}$.
We deduce the compatibility properties required by applying
to $\pi_ { F(\mathcal{Q})}$ :

- if $\lambda$ is not above $q$, the theorem of Carayol
(\cite{Carayol}) completed by Taylor (\cite{Tay});

- if we are in the case (i)  and $\lambda$ is above $q\not= 2$,
the theorems of Breuil (\cite{Breuil-cong}) and Berger
(\cite{Berger}) to get that $\rho$ is crystalline at $q$;

- if we are in the case (ii), the theorem of Saito (\cite{Saito})
because we know that $\pi_ { F(\mathcal{Q})}$ is Steinberg
at one prime of $F(\mathcal{Q})$
(note that if $\rho$ is unramified outside $p$, $\rho$ is semistable
not crystalline of  weight 1 at $p$, and it follows from
\cite{Breuil-cong} and  \cite{Berger} that  $\pi_ { F(\mathcal{Q})}$
is not unramified at primes above $p$).

The last statement
of (ii) follows from the arguments used for Corollary E, or Corollary
2.4, of
\cite{[T02]} which use results of Blasius-Rogawski (\cite{BR}). Let us do it when
$k(\rhobar)=2$ or $p+1$ (and we consider weight 2 semistable lifting
$\rhobar$). By \cite{BR}, we know that the restriction of $\rho$
to $G_F$ is a direct factor of the Tate module of an abelian variety
defined over $F$. By Galois descent, $\rho$ also is a direct factor
of the Tate module $V_p (A)$ of an ablian variety $A$ over $\Q$,
which we may assume is simple. Let $E_A$ be the center of the
skew field $\mathrm{End}_{\Q} (A)$.
Replacing $A$ by a simple factor of
 $A\otimes_{E_A} E'$ for $E'$ a finite extension
of $E_A$ that splits $\mathrm{End}_{\Q} (A)$, we may suppose that
 $E_A =\mathrm{End}_{\Q} (A)$.
 By Faltings, we have an embedding of $i$ of
$E_A$ in $\aQp$ such that $V_p (A)\otimes _{E_A ,i} \aQp$ is isomorphic
to $\rho$. Enlarging $E$ so it contains $E_A$ and  replacing
$A$ by $A\otimes_{E_A} E$, we get our abelian variety. It is
semistable over $\Q$ by the properties of $(\rho _{\lambda })$.

\end{proof}

\noindent{\bf Remark:}
It is easy to see that all members of the compatible system are
irreducible.
We say that the compatible system has good reduction at a place $v$ if
the Weil-Deligne parameter at $v$, $(\tau,N)$, is such that $\tau$ is
unramified and $N=0$, and we say it has semistable reduction at $v$ if
$\tau$ is
unramified. We will say that the $(\rho_{\lambda})$ of the theorem
give a lift of $\rhob$ to a {\it minimal strictly compatible system}.

\section{Low levels and weights}

Theorem \ref{minimal} when combined with modularity lifting results in
\cite{Wiles}, \cite{SW2}, and the theorems of Fontaine,
\cite{Fontaine},
together with  their generalisations due to Brumer and Kramer, and Schoof, \cite{BK},
\cite{Schoof1},  and \cite{Schoof2},
has a number of corollaries.

Part (i) of the following consequence of his conjectures was spelled out by
Serre in Section 4.5 of \cite{Serre3} and which we now prove unconditionally.

\begin{theorem}\label{lowlevel}
 (i) There is no finite flat group scheme over $\Z$ of $(p,p)$ type which is
 irreducible. (In fact there is no irreducible, odd 2-dimensional
 representation
 $\rhob$ that is unramified outside $p$ and whose Serre weight
 $k(\rhob)$ is 2). 

(ii) Assume $p>3$. There is no semistable $\rhob$ with
 (prime to $p$) Artin conductor either 2,3,5,7,13 
 and $k(\rhob)=2$.

 (iii) For $p=5,7,13$ there is no Serre type $\rhob$ that is
 unramified outside $p$ and such that $k(\rhob)=p+1$.

\end{theorem}

\begin{proof} Let us prove (i). 
The case $p=3$ is taken care of by \cite{Serre2}, page 710, which proves
that
there is no odd 2-dimensional irreducible representations of
$\Gal$ in $\mathrm{GL}_2 ({\overline{\mathbb{F}_{3}}})$
unramified
outside $3$.
For $p>3$,  (i)
follows from Theorem \ref{minimal} and the main result
of \cite{D} and \cite{Wint}.
For instance \cite{Wint} uses a 7-adic representation
that is part of the system to apply Fontaine's result, while \cite{D} uses
a 3-adic  representation
that is part of the system and then uses
\cite{SW2} and Serre's result from \cite{Serre2}.

For (ii), the argument is similar. This time again
use
Theorem
\ref{minimal} to get minimal compatible lift $(\rho_\lambda)$ of
$\rhob$ (of weight 2), and
use the refinements of \cite{Fontaine}
in Brumer-Kramer and Schoof, \cite{BK} and \cite{Schoof1}, which yield that
there are
no semistable abelian varieties over $\Q$ with good reduction
outside one of   2,3,5,7 or 13.

Part (iii) again follows by quoting results in loc. cit. after using
Theorem \ref{minimal} to get weight 2 semistable liftings.

\end{proof}

\noindent{\bf Remark:} Assuming GRH and assuming $p>3$,
one deduces using Calegari \cite{Calegari} that there is no irreducible
$\rhob$
with
conductor 6 and Serre weight 2. By similar methods, using
Theorem 1.3 of \cite{Schoof1} one can probably rule out some more
cases.
Theorem
\ref{lowlevel}  gives a slightly simpler proof of Fermat's Last Theorem
which
was proven in \cite{Wiles} (we will not have to use the most difficult of
the
level lowering results which are
those proven by Ribet in \cite{Ribet} (at least for $p>49$ to avoid
residually solvable cases), nor the
modularity
of the putative $\rhob$ of level 2, and weight 2 which a solution of
Fermat's equation leads to!)

\begin{cor}\label{groupschemes}
 If $p$ is odd then the only $(p,p)$-type
 finite flat group schemes over $\Z$ are $\Z/p\Z \oplus \Z/p\Z$,
 $\Z/p\Z \oplus \mu_p$ or $\mu_p \oplus \mu_p$.
\end{cor}

For $p=2$, see Abrashkin \cite{[A89]}.

\begin{proof}
After part (i) of Theorem \ref{lowlevel} this follows from Serre's
arguments in Section 4.5 of
\cite{Serre3}. 
\end{proof}

In fact using our methods we can also rule out the existence of some higher
weight $\rhob$ in accordance with the predictions of Serre.

\begin{theorem}\label{lowweight}
 There is no Serre type $\rhob$ of level 1
such that $2 \leq k(\rhob) \leq 8$, or
 $k(\rhob)= 14$.

\end{theorem}

\begin{proof} We begin by observing that the Serre weight of $\rhob$  is
even. We begin by noting that there is no irreducible
$\rhob$ with $p=3$ and unramified outside 3. This in fact has been proven
directly By Serre in \cite{Serre2}, page 710,
we know   that there is no irreducible
$\rhob$ with $p=3$ and unramified outside 3.

\begin{itemize}

\item $k(\rhob)=4$:
Now let $p$ be any prime $>3$, and $\rhob$ be an irreducible
representation with $k(\rhob)=4$.
We use Theorem \ref{minimal} to get a strictly compatible system
$(\rho_\lambda)$, of weight 4 (with good reduction everywhere).
Consider a prime $\lambda$ above 3: then by
results in 4.1 of \cite{BLZ}, and as the residual
representation is
globally and hence of course locally at 3 reducible, we see that
$\rho_\lambda$ is ordinary at 3, and hence by Skinner-Wiles \cite{SW2} (note
that
the 3-adic representation we are considering is 3-distinguished) corresponds
to a cusp form of level 1 and weight 4 of which there are none.

\item $k(\rhob)=6$: Suppose we have an irreducible $\rhob$ with
$k(\rhob)=6$ as in the theorem ($p>3$).

We begin by proving the claim that there is no irreducible
$\rhob$ with $p=5$ and unramified outside 5 (this was known earlier only
under
the GRH in \cite{Brue}). By twisting we may assume
that $k(\rhobar) \leq 6$, and thus by 
     the previous step and theorem \ref{lowlevel} that $k(\rhobar)=6$.
Using Theorem
\ref{minimal}, we get a strictly compatible sytem of weight 2
representations $(\rho_\lambda)$ (whose Weil-Deligne parameters are
unramified
outside 5 and semistable at 5)
such that for a place above 5 the residual representation is
$\rhob$. But again by what was recalled in proof of
Theorem \ref{lowlevel}, such a system cannot exist by results of
\cite{Schoof1}.
We thus see that we have proved the claim.

Now let $p>5$ be any prime, and $\rhob$ be an absolutely irreducible 
representation with $k(\rhob)=6$. This time we use Theorem \ref{minimal} to get a
compatible system $(\rho_\lambda)$ of weight 6, i.e., Hodge-Tate of
weights
$(0,5)$, and with good reduction everywhere.
Consider a prime $\lambda$ above 5: then by results in \cite{BLZ},
and as the residual representation is globally and hence of course locally
reducible at 5, we see that
$\rho_\lambda$ is ordinary at 5, and hence by Skinner-Wiles \cite{SW2}
corresponds
to a cusp form of level 1 and weight 6 of which there are none.

\item $k(\rhob)=8$: Consider an irreducible $\rhob$ with $k(\rhob)=8$
  as in the theorem.
By what we just did we see that $p>5$.

We begin by proving the claim that there is no such irreducible
$\rhob$ with $p=7$ and unramified outside 7.
By twisting we may assume
that $k(\rhobar) \leq 8$, and thus by the previous steps that $k(\rhobar)=8$. 
Using Theorem
\ref{minimal}, we get a strictly compatible sytem of weight 2
semistable representations $(\rho_\lambda)$ (with good reduction
outside 7, and semistable reduction at 7)
such that for a place above 7 the residual representation is
$\rhob$. But again by what was recalled in proof of
Theorem \ref{lowlevel}, results of \cite{Schoof1} show that
such a system cannot exist.
We thus see that we have proved the claim.

Now let $p$ be a prime $>7$, and $\rhob$ be an irreducible
representation with $k(\rhob)=8$. This time we use Theorem
\ref{minimal} to get a compatible system $(\rho_\lambda)$ of weight 8, i.e.,
of
Hodge-Tate weights
$(0,7)$ (with good reduction everywhere).
Consider a prime $\lambda$ above 7: then by results in \cite{BLZ},
and as the residual representation is globally and hence of course locally
reducible at 7, we see that
$\rho_\lambda$ is ordinary at 7, and hence by Skinner-Wiles corresponds
to a cusp form of level 1 and weight 8 of which there are none.

\item $k(\rhob)=14$: Consider an irreducible
$\rhob$ with $k(\rhob)=14$ as in the theorem.
We first deal with the case $p=11$. It is easy to see, by the definition of
the
Serre weight, that $\rhob$ is
a twist by $\chi$ of a level 1 representation with Serre weight 2,  that as
we
have ruled out.

We next prove the inexistence of $\rhob$, with
$k(\rhob)=14$,
for $p=13$.
Using Theorem \ref{minimal}
we get a strictly compatible sytem of weight 2
semistable representations $(\rho_\lambda)$ (with good reduction
outside 13 and semistable at 13)
such that for a place above 13 the residual representation is
$\rhob$. But again by what was recalled in proof of
Theorem \ref{lowlevel}, such a system cannot exist.
We thus see that we have proved the inexistence of  $\rhob$.

Now let $p > 11$, and $\rhob$ be a mod $p$ Serre type
representation of level 1, with $k(\rhob)=14$. This time we use Theorem
\ref{minimal} to get a compatible system $(\rho_\lambda)$ that has good
reduction everywhere and is crystalline of weight 14, i.e., of
Hodge-Tate weights $(0,13)$.
Consider a prime $\lambda$ above 13. Then by results in \cite{BLZ},
and as the residual representation is globally and hence of course locally
reducible at 13, we see that
$\rho_\lambda$ is ordinary at 13, and hence by Skinner-Wiles \cite{SW2}
corresponds
to a cusp form of level 1 and weight 14 of which there are none.

\end{itemize}

\end{proof}

The following corollary we have proved in the course of the proof above is
worth
noting (the cases $p=2,3$ are theorems of Tate and Serre):

\begin{cor}\label{small}
For the primes $p=2,3,5,7$ there are no irreducible, odd
2 dimensional mod $p$ representations of $\Gal$, unramified outside $p$.

\end{cor}

\begin{theorem}\label{lowweight1}

 For a odd prime $p$, there is no 2 dimensional odd irreducible $p$-adic
representation of $\Gal$ which is unramified outside $p$ and crystalline
 of Hodge-Tate
weights $(0,w)$, with $0 \leq w \leq 7$ or $w=13$ and $p \geq w$
(i.e., we rule out
$w=1,3,5,7,13$ as $w$ is easily seen to be odd \cite{Wint}). 
There is no irreducible strictly compatible
system
of 2 dimensional irreducible $\lambda$-adic representations of $\Gal$
which has good reduction everywhere and is crystalline
of Hodge-Tate
weights $(0,w)$, with $0 \leq w \leq 7$ or $w=13$.

\end{theorem}

\begin{proof} Let $k=w+1$.
We prove the statement that
there is no 2 dimensional odd irreducible $p$-adic representation $\rho$
which is unramified outside $p$ and crystalline  at  $p$ of Hodge-Tate
weights $(0,w)$ with $1 \leq w \leq 7$ or $w=13$, and in each case $p \geq
w$,
from which the second statement follows easily.
If $w=1$, or $w=3$ and $p\geq 7$, this is already in \cite{D} and \cite{Wint}.
Notice that in each
of the other cases, for the values of $w$ considered, $w$ is an odd prime.
The case $p=3$ has already been taken care of 
by \cite{Serre2}, \cite{SW2} and \cite{BLZ}
(its relevant only for
$k=1,3$).

If the representation is residually reducible, then again using \cite{BLZ}
(which can
be used because of our assumption $p+1\geq k$)
and \cite{SW2} we see that $\rho$ arises from a  cusp form of level
1 and weight $k$ either
at most 10 or weight 14, which do not exist. If residually the
representation is
irreducible,
we get a compatible system as in Theorem \ref{minimal}, and then
consider the corresponding residual representation
$\rhob$ at the prime $w$. We claim that $\rhob$ is reducible. 
We know that locally
at $w$ by \cite{BLZ} that the Serre
weight of the residual representation
can be either 2 or $w+1$ both of which we have ruled out
as weights which can occur for irreducible
$\rhob$. Hence the residual representation is reducible
and we are done by the previous analysis.
\end{proof}

\begin{theorem}\label{weight12}
Let  $\rhob$ be a Serre type, level 1 representation,
with Serre weight 12. Then $\rhob$ arises from the Ramanujan $\Delta$
function.
\end{theorem}

\begin{proof}
From Corollary \ref{small} we see that we may assume
that $p \geq 11$. Using Theorem \ref{minimal},
we get a strictly compatible system $(\rho_\lambda)$
with good reduction everywhere and is
crystalline of Hodge-Tate weight $(0,11)$ and at a prime above $p$ the
representation is residually $\rhob$. Consider a prime above
11 and reduce the corresponding 11-adic representation that is a part of
$(\rho_\lambda)$ mod 11. If we get a reducible representation then again
using
Remark 4.1.2 of \cite{BLZ} and Skinner-Wiles \cite{SW2} we are done.
Otherwise
the
representation,
call it $\rhob_{11}$ is irreducible and again using \cite{BLZ}
its easily seen  as before to
be tr\`es ramifi\'ee (and hence of Serre weight 12), as we have ruled
out the weight 2 case.

Now using Theorem
\ref{existminimallift}, we
get another lift $\rho '_{11}$ of $\rhob_{11}$, which is
unramified outside 11 and semistable at 11.
By Theorem \ref{minimal} $\rho '_{11}$ arises
from an abelian variety $A$ defined
over $\Q$. Then, as
$\rho '_{11}$ is unramified outside 11 and semistable at 11,
$A$ is semistable and has good reduction outside 11.
But Schoof \cite{Schoof1} has proven that such an abelian
variety $A$ is isogenous to a power of $J_0(11 )$.
The Galois representation on points of order 11 of
$J_0 (11 )$ is absolutely irreducible, is ordinary at 11, and is isomorphic
to the representation modulo 11 associated to $\Delta$
(for example Th\'eor\`eme 11 of \cite{Serre2} for the latter), and
thus $\rhob_{11}$ itself arises from $\Delta$.
The image of the mod 11 representation arising from $\Delta$, and
hence that of $\rhob_{11}$, is all of $GL_2(\F_{11})$
(see \cite{Serre4}).
Thus using modularity theorems of Wiles \cite{Wiles} for the 11-adic
representation that figures in the compatible system $(\rho_\lambda)$ (and
we
also need as before to use \cite{BLZ}, to see that the lift is ordinary at
11
as the residual weight is 12,
and hence $\rhob_{11}$ is tr\`es ramifi\'ee and in particular ordinary),
we conclude that the compatible system $(\rho_\lambda)$ arises from the
Ramanujan $\Delta$ function, and hence that $\rhob$ arises from the
Ramanujan
$\Delta$ function.
(Of course sometimes there may be no irreducible
$\rhobar$ as in the theorem: but these primes
have been determined by Swinnerton-Dyer, they being 2, 3, 5, 7, 691,
see \cite{Serre4}.)
\end{proof}

\begin{cor}
Any 2 dimensional irreducible $p$-adic
representation $\rho$ which is unramified outside $p$ and crystalline
at  $p$ of Hodge-Tate
weights $(0,11)$, and $p \geq 11$, arises
from the Ramanujan $\Delta$ function. Furthermore any
irreducible strictly compatible system
of 2 dimensional irreducible $\lambda$-adic representations of $\Gal$
with good reduction everywhere and of Hodge-Tate
weights $(0,11)$ arises from the Ramanujan $\Delta$ function.
\end{cor}

\begin{proof} We prove only the first statement. If $\rho$ is residually
reducible, using
\cite{BLZ} and  \cite{SW2} we are done.
Assume not. If $p\geq 13$, as $2\times 12\not=p+3$,
we can apply theorem 6.1. of \cite{[T03]}, and by \cite{D}
and \cite{Wint} we can make $\rho$
of a compatible system $(\rho_{\lambda})$ and consider a representation
above 11 of this system. 
We are reduced to the case $p=11$. 
If residually it is reducible
we see as before using \cite{SW2} that it, and hence the compatible
system, arises from a newform, and we are done. If it is irreducible,
then by  Theorem \ref{weight12} we know that the residual
representation $\rhob$ arises from the $\Delta$ function, and we are
done by \cite{SW3}.
\end{proof}

\begin{cor}\label{ugly}
For $p=11,13$ the only $\rhob$ of level 1 of Serre type are
those that have a twist of Serre weight 12  and these arise from
$\Delta$ (up to twist), and  those $\rhob$
one of whose twists is weight 10,
and such that the local representation at $p$ is non-semisimple.

For $p$ bigger than 13, there is no Serre type $\rhobar$ of level 1 of
Serre weights
$p-1,p-3,p-5,p-7,p-13$ when the local  representation at $p$ is
reducible
and split and in the case $p-11$ if it exists it arises from $\Delta$
up to twist.
There is no Serre type of $\rhobar$ of weights
$p-1,p-3,p-5,p-11$ when the local  representation at $p$ is
irreducible
and in the case $p-9$ if it exists it arises from $\Delta$ up to twist.
\end{cor}

\begin{proof} This follows from the previous work using the
observation that when  $\rhob|_{D_p}$ is semisimple with $2 \leq k(\rhob) <
p$, $\rhob \otimes \overline{\chi_p}^i$  is such that $2 \leq k(\rhob \otimes
\overline{\chi_p}^i) \leq p+1$ for some $i$ with $0 < i <p-1$.
When the the local  representation at $p$ is
reducible and split, such twist is of Serre weight $p-k(\rhob )+1$.
When the local  representation at $p$ is 
irreducible and $k(\rhob )\not=2$, such a twist is of weight $p-k(\rhob )+3$.
For
example, for $p=11$ in the locally irreducible case, we use our result for
Serre
weight 4,
and in the locally completely reducible case for Serre weight 2. For $p=13$
in the irreducible case we use our result for Serre weight 6, and in the
completely reducible case for Serre weight 4.
\end{proof}

\noindent{\bf Remark:} The case of Serre type $\rhob$ of level 1 with
$k(\rhob)=10$, which should not exist, cannot be treated directly by
the methods here (as 9 is not a prime for
instance). In this case it might be possible to get potentially
Barsotti-Tate lifts at $p$ using methods of this paper and
\cite{Conrad}, and then using results of \cite{Schoof2} which are
conditional on GRH that rule out abelian varieties that have
everywhere
good reduction over $\Q(\zeta_{11})$.

\section{Two reductions of residual modularity to modularity lifting
results,
via induction on primes}

\subsection{Level 1 case of Serre's conjecture: reduction
to modularity lifting for moderate weights}

We now address the question of proving Serre's conjecture for Serre
type
$\rhob$ of level 1. By results of Tate and Serre we may assume $p>3$.
By twisting we may assume that the Serre weight $k(\rhob) \leq p+1$. Using
Theorem \ref{minimal} we get a strictly compatible system
$(\rho_\lambda)$ with good reduction everywhere and crystalline of weight
$k(\rhobar)$ and such that a place above $p$ the residual representation is
$\rhob$. Now if we reduce
the compatible system at a prime above 3 by Serre's result quote before
we get a residually reducible representation and
hence one might try to adapt \cite{SW2} to this per force non-ordinary and
very
high
weight situation! This might be very hard technically.
We ameliorate this a little by proving the following theorem which reduces
Serre's conjecture in the level 1 situation to what might hopefully be a
more tractable
modularity lifting theorem. The argument of moving the prime around was
suggested by the proofs in the section
on low levels and low weights (especially proof of Theorem \ref{lowweight}).

In both the theorems below because of the condition of oddness and the
semistability assumption the Serre weight of the representations
considered is always even and in particular $\neq p$.

\begin{theorem}\label{bertrand}
 For every odd prime $p$, assume the following statement:

Let $\rhob:\Gal \rightarrow GL_2(\FF)$ be a continuous representation of
level
one with $\FF$ a finite field of characteristic
 $p$, that is either irreducible and modular, or is reducible.
 For any $\rho:\Gal \rightarrow GL_2({\cal O})$, with ${\cal O}$ the ring of
integers of a finite extension of $\Q_p$,
 that is a continuous lift of such a $\rhob$, and
 that is unramified outside $p$ and crystalline
 of weight between 2 and $2p$ (in fact may even assume that the weight is
actually at most 1 more than the next prime after $p$), assume $\rho$ is
modular.

Then Serre's conjecture in level 1 is true.
\end{theorem}

\begin{proof}

We prove this theorem by induction on the prime $p$. (In fact what the
argument
will prove is that if we know the lifting
statement for a prime, and Serre's conjecture for that prime, then we know
the
level one Serre conjecture for the next prime, or even a much larger batch
of
subsequent prime(s).)
As we have proven Serre's conjecture in the level 1 case for primes less
than
11, we start with 7 for which by Corollary \ref{small} we know
the level 1 case of Serre's conjecture. Suppose
Serre's level one conjecture is proven for a prime $p_n$. We want to prove
it
for the next prime
$p_{n+1}$. Thus assume we have an irreducible, odd, 2-dimensional mod
$p_{n+1}$
representation $\rhob$ of $\Gal$ which by twisting we can assume
has Serre weight $k(\rhob) \leq p_{n+1}+1$.

Now we use Theorem \ref{minimal} to get a compatible system
$(\rho_\lambda)$
that is unramified everywhere, is crystalline of weight $k(\rhob)$ and for a
prime above $p_{n+1}$  we get a $p_{n+1}$-adic representation
that is part of the compatible system, and which lifts $\rhob$. This $\rho$
is
unramified outside $p_{n+1}$ and is crystalline at $p_{n+1}$ of Hodge Tate
weights $(0,k(\rhob)-1)$. We use Bertrand's postulate to see that $p_{n+1}
\leq
2p_n-1$, and thus by assumption we get the modularity
of the member of the compatible system $(\rho_\lambda)$ at a prime above
$p_n$
by the hypothesis of the theorem, as the induction hypothesis guarantees
residual modularity for mod $p_n$ representations. Thus we get the
modularity of
the compatible sytem $(\rho_\lambda)$ and thus that of
$\rhob$. This completes the induction step and hence the proof of the
theorem.
\end{proof}

\noindent{\bf Remarks:}

- The reduction to moderate weights (i.e., between 2 and $2p$)
in Theorem \ref{bertrand}
perhaps is technically critical as the conjectures of Breuil
in \cite{Breuil} (see Conjecture 6.1)
about reductions of crystalline representations of weights at most $2p$ have
a
simple form. These
conjectures may be close to be proven, and this might be helpful in proving
the
hypothesis of Theorem \ref{bertrand}.

- When the liftings $\rho$ are ordinary the hypothesis of Theorem
\ref{bertrand} is satisfied by results of \cite{Wiles}, \cite{SW2} (note
that as
weights are
even
the distinguishedness hypothesis of that paper is satisfied).

\subsection{Killing ramification}

The process of killing ramification is the following.
Suppose you wish to prove that a compatible system
$(\rho _{\lambda })$ is modular. Let $\lambda_0$ be above
a prime of ramification of $(\rho _{\lambda })$. One 
applies the theorem  \ref{minimal} to a cyclotomic twist of
$\rhob _{\lambda_0 }$ to get a compatible system 
$(\rho ' _{\lambda ' })$ whose set of ramification
primes is smaller than the set of ramification primes
of $(\rho _{\lambda })$. If one knows by induction modularity of
$(\rho ' _{\lambda ' })$, one get modularity of 
$\rhob _{\lambda_0 }$, hence modularity of 
$(\rho _{\lambda })$ if one has the needed modularity lifting theorem.
We give an example of a more precise statement:

\begin{theorem}\label{levelone}

For every odd prime $p$ we assume the following statement:

Let $\rhob':\Gal \rightarrow GL_2(\FF)$ be a continuous representation of
squarefree conductor with $\FF$ a finite field of characteristic
 $p$, that is either irreducible and modular, or is reducible.
 For any $\rho':\Gal \rightarrow GL_2({\cal O})$, with ${\cal O}$ the ring
of
integers of a finite extension of $\Q_p$,
 that is a continuous lift of such a $\rhob'$,  and which is semistable at
primes outside $p$ and of Serre
weight between 2 and $p+1$ at $p$, assume $\rho'$ is modular.

From this we deduce that Serre's conjecture for level one implies that
any $\rhobar$ (with $p$ bigger than 3) which
is an irreducible 2-dimensional odd
mod $p$ representation with odd squarefree conductor prime to 3,
and whose determinant is ramified only at $p$, and such that its Serre
weight
$k(\rhobar)$ is either $p+1$, or is at most 1 more than the least prime
ramified
in $\rhobar$, is modular.
\end{theorem}

\begin{proof}
Consider a semistable $\rhob$ which is an irreducible
mod $p$ representation with $p$ an odd prime. (In the process below
whenever we reach a representation that is residually solvable we can stop.)

- We first begin by showing
how we may deal with the case $k(\rhob)=p+1$, and reduce it to the other
allowed
cases. We use
Theorem \ref{minimal} to lift $\rhob$ to a strictly compatible
semistable system
$(\rho_\lambda)$ of weight 2 that is has semistable reduction at  $p$. Consider a large prime $q$ at which a $q$-adic
member of
the system is crystalline of Hodge-Tate weights $(1,0)$ and the residual
representation $\rhob_q$ is irreducible. It will
be
enough
to prove modularity
of $\rhob_q$ as then known modularity lifting theorems as in \cite{Wiles}
would prove modularity of the entire compatible system $(\rho_\lambda)$ and
hence of $\rhob$. 

- Let $S=\{q_1,\cdots,q_n\}$
be the primes that are ramified for $\rhob$ not including $p$, written
in increasing order, and we assume
that  $k(\rhob)$ is at most 1 more than the least prime ramified in
$\rhobar$.

The proof is by induction on the cardinality of $S$ for the type of
$\rhobar$
in the statement. The case when $S$ is empty is the case
we are assuming.
Use  Theorem \ref{minimal} to get a strictly compatible minimal
system that we again denote by $(\rho_\lambda)$ which at a prime above
$p$ reduces to $\rhob$, and is crystalline of weight $k(\rhob)$. We consider
a
residual representation $\rhob_{Q_1}$
attached
to this system at a prime $Q_1$ above $q_1$.
The Serre weight of
$\rhob_{Q_1}$ is $\leq q_1+1$ by assumption and the (prime to $q_1$)
Artin conductor of $\rhob_{q_1}$ is divisible by at least one prime
less than that of the (prime to $p$)
Artin conductor of $\rhob$. Then by the inductive hypothesis,
we deduce that $\rhob_{Q_1}$ is modular, and then by the lifting
hypothesis of the theorem we
see that $\rho_{Q_1}$, and hence $(\rho_\lambda)$, arises from
a newform.

- In the end we are  reduced, assuming the hypothesis of the theorem, to
proving the modularity of $\rhob$ when it is
irreducible of level 1 and $p$ is some odd prime as asserted.

\end{proof}

\noindent{\bf Remarks:}

- The modularity lifting hypothesis of Theorem \ref{levelone}
 maybe accessible at least for weights at most $p-1$ using results of
\cite{BM} (and that explains why we have made the restrictive weight
 hypothesis in the theorem).

- To us the cases of modularity lifting that are needed  in Theorem
\ref{bertrand}
and Theorem \ref{levelone} that seem hardest, are the residually degenerate
cases, i.e., dihedral induced from $\Q(\sqrt { {(-1)} ^{p-1/2} p})$ or
reducible and when the lifts whose modularity needs
to be established are non-ordinary liftings. The residually
``non-degenerate''
cases, while not as yet proven or available in the literature, seem
accessible
because of the basic method of Wiles et al. and its recent
developments due to Kisin
in \cite{Kisin}, together with the results
of Breuil, Berger, Li, Mezard, Zhu, see
\cite{Breuil}, \cite{Berger}, \cite{BLZ}, \cite{BM}.

\nocite{*}
\bibliographystyle{plain}

\noindent{\bf Addresses of the authors:}

\noindent CK:
155 S 1400 E,
Department of Mathematics,
University of Utah,
Salt Lake City,
UT 84112,
USA. e-mail: {\tt shekhar@math.utah.edu}, {\it and},
\noindent School of Mathematics,
TIFR,
Homi Bhabha Road,
Mumbai 400 005,
INDIA. e-mail: {\tt shekhar@math.tifr.res.in}

\vspace{3mm}

\noindent JPW:
Universit\'e Louis Pasteur
D\'epartement de Math\'ematiques, IRMA
7, rue Ren\'e Descartes
67084 Strasbourg Cedex
France. e-mail: {\tt wintenb@math.u-strasbg.fr},
tel 03 90 24 02 17 , fax 03 90 24 03 28

\end{document}